\documentclass[letterpaper, 10 pt, conference]{ieeeconf}  

\IEEEoverridecommandlockouts                              
\usepackage{amsmath,amssymb,amsfonts}

\usepackage{amsthm}
\usepackage{setspace}
\usepackage{xcolor}
\usepackage{graphicx}   
\usepackage[caption=false]{subfig}

\newtheorem{thm}{Theorem}

\newtheorem{cor}{Corollary}
\newtheorem{lem}{Lemma}

\newtheorem{assum}{Assumption}

\newtheorem{defn}{Definition}
\newtheorem{exmp}{Example}
\newtheorem{rem}{Remark}
\newtheorem{prob}{Problem}

\newcommand{\Real}{\mathbb{R}}

\newcommand{\setC}{\mathcal{C}}
\newcommand{\setD}{\mathcal{D}}
\newcommand{\boundC}{\partial\mathcal{C}}
\newcommand{\intC}{\text{Int}(\mathcal{C})}

\newcommand{\classKinf}{\mathcal{K}_\infty}

\DeclareMathOperator*{\argmin}{\arg\min}

\newcommand{\fdsb}[1]{\textcolor{black}{#1}}

\title{\LARGE \bf
Provably Safe Control of Lagrangian Systems \\ in Obstacle-Scattered Environments
}

\author{Fernando S. Barbosa$^{1}$, Lars Lindemann$^{1}$, Dimos V. Dimarogonas$^{1}$ and Jana Tumova$^{1}$
\thanks{*Work supported by EU H2020 Co4Robots, by the Swedish Research Council (VR) and by the Knut och Alice Wallenberg Foundation (KAW).}
\thanks{$^{1}$KTH Royal Institute of Technology, Stockholm, Sweden.
        {\tt\small \{fdsb, lindem, dimos, tumova\}@kth.se}}%
}

\begin{document}

\maketitle
\thispagestyle{empty}
\pagestyle{empty}

\begin{abstract}
We propose a hybrid feedback control law that guarantees both safety and asymptotic stability for a class of Lagrangian systems in environments with obstacles. 
Rather than performing trajectory planning and implementing a trajectory-tracking feedback control law, our approach requires a sequence of locations in the environment (a path plan) and an abstraction of the obstacle-free space.
The problem of following a path plan is then interpreted as a sequence of reach-avoid problems: the system is required to consecutively reach each location of the path plan while staying within safe regions. 
Obstacle-free ellipsoids are used as a way of defining such safe regions, each of which encloses two consecutive locations.
Feasible Control Barrier Functions (CBFs) are created directly from geometric constraints, the ellipsoids, ensuring forward-invariance, and therefore safety. 
Reachability to each location is guaranteed by asymptotically stabilizing Control Lyapunov Functions (CLFs).
Both CBFs and CLFs are then encoded into quadratic programs (QPs) without the need of relaxation variables.
Furthermore, we also propose a switching mechanism that guarantees the control law is correct and well-defined even when transitioning between QPs. Simulations show the effectiveness of the proposed approach in two complex scenarios.
\end{abstract}

\section{INTRODUCTION}

The task of autonomous robot navigation in complex environments requires several layers of planning and control to work together and in harmony. 
Traditionally, path planning concerns longer-term planning in (discrete) high-level representations of the robot's environment. Motion planning is a mid-term layer, in charge of generating (optimal) trajectories in the state-space of the robot in accordance with a path plan. Feedback control deals with short-term navigation, typically tracking the (optimal) trajectory and handling only local obstacle avoidance, but leaving the challenges of the structured environment to the higher planning layers. Integrating planning and control for safety-critical systems is important for guaranteeing safety is kept even among disturbances and unmodeled dynamics, therefore ensuring safety in both short and long term.


\begin{figure}
    \centering
    \includegraphics[width=0.85\linewidth]{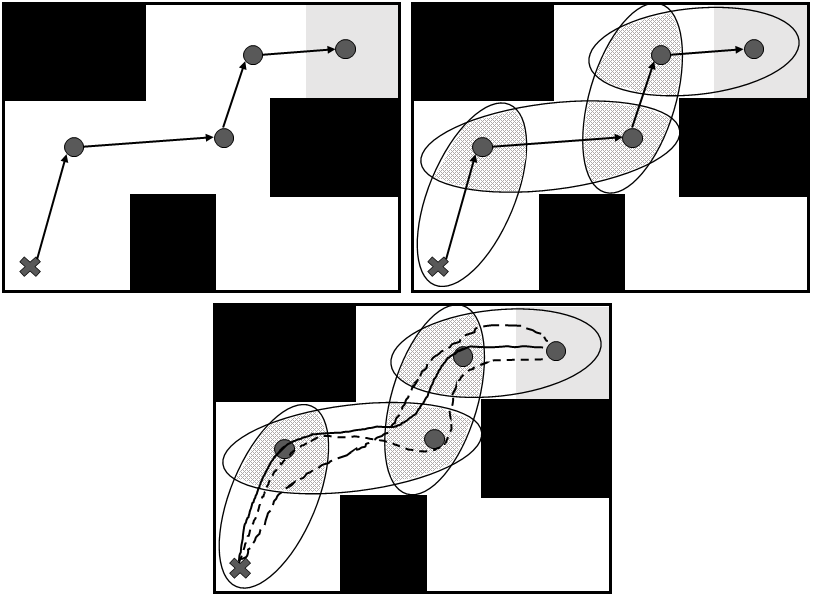}
    \vspace{-3mm}
    \caption{Illustration of the proposed approach. The top left figure depicts a 2D workspace with obstacles (black rectangles), a goal region (gray), an initial position of a Lagrangian system (marked with `x'), and a path plan (circles). The top right figure shows \emph{safe sets} (obstacle-free ellipsoids) enclosing each edge, highlighting the intersections of two consecutive ones. The bottom figure shows three trajectories a system can develop depending on its dynamics and parameters.}
    \label{fig:example}
\end{figure}


Safety requires undesired states to never be reached, and can be seen as the dual problem of reachability, which requires desired states to be eventually reached \cite{ames2017control,romdlony2016stabilization,ames2019control}. The main goal of this work is to propose a method to ensure safety and asymptotic stability of Lagrangian systems in order to fulfill a path plan in environments containing obstacles; a path plan is given as a sequence of locations (positions) in the configuration space.

Several related papers treated the theme of safe motion planning and control of dynamical systems lately. In \cite{liu2017planning}, the problem is posed to plan dynamically feasible trajectories guaranteed to be inside the obstacle-free space, but with the rather restrictive assumption that a low-level controller capable of tracking such trajectory exists. 
In \cite{herbert2017fastrack}, Hamilton-Jacobi reachability analysis of a pursuit-evasion game between the system dynamics and its simplified version is used in order to find tracking error bounds, which are then used to generate safety bubbles around the planned trajectory. 
Since such an approach can be overly conservative, it is extended in \cite{fridovich2018planning} to allow planning for two situations: high velocity with large bounds, or slow with small bounds. 
Lastly, contraction theory and convex optimization are used in \cite{singh2017robust} to calculate a fixed-size tube around a nominal trajectory within which the system is guaranteed to remain; however, it assumes the existence of a kinodynamic planner (such as \cite{lavalle2001randomized}) and builds on robust trajectory tracking formulation. Although these works are the closest to ours, we neither pose our problem as trajectory-tracking nor do we require (possibly computationally expensive) offline steps in order to calculate error bounds. Instead, our approach offers a feedback control law that drives the system through the obstacle-free subspace of the environment.

Polytopic-trees \cite{sadraddini2018sampling} and LQR-trees \cite{tedrake2010lqr} are feedback motion planners that aim at building stabilized sets on the state-space of the system that cover the entire state-space, ensuring convergence of the system to the goal region. Both approaches grow trees backwards, rooted at the goal region, thus not applicable to receding-horizon problems. Furthermore, their applicability is limited by the sums-of-squares verification, which scales poorly with dimension.

On the feedback control theory side, Control Barrier Functions (CBFs) were proposed as a manner to provide and ensure safety posed as set invariance of a system.
CBFs, alongside their dual Control Lyapunov Functions (CLFs), have been applied to adaptive cruise control \cite{ames2014control}, systems evolving on manifolds \cite{wu2015safety}, multirobot systems \cite{wang2017safety}, fixed-wing aircraft \cite{squires2018constructive}, and others. 
The works in \cite{nilsson2018barrier,srinivasan2018control,garg2019control} use CBFs for planning and control by considering general system dynamics and assuming the CBF to be valid.  These works further pose the conjunction of control barrier functions, complicating the synthesis of valid CBFs in practice. Indeed, it turns out that finding valid CBFs is, in general, a difficult problem. A related problem is to find barrier certificates for dynamical systems. Sum-of-squares approaches for polynomial systems, hence not necessarily including the class of Lagrangian systems, have appeared in this respect in \cite{prajna2007framework}. We, on the other hand, use specific properties of the system at hand and select, in a useful way, convex barrier functions showing that the CBF is feasible, ensuring the correctness of our approach.

Our proposed approach is illustrated in Fig.~\ref{fig:example}. Given a path plan, we propose a way of encoding safety as CBFs directly from geometric constraints, and asymptotic stability as CLFs. More precisely, we propose a solution that interprets the task as a sequence of reach-avoid problems, where we explore the use of (obstacle-free) ellipsoids as safe-sets enclosing two consecutive configurations. A closed-loop controller is obtained by solving a Quadratic Program (QP) encoding both CBF and CLF constraints without the need for a relaxation variable. Furthermore, a switching mechanism is proposed for transitioning, in a correct manner, from one safe-set to another, so that the obtained controller is well-defined even during these transitions, and therefore fulfilling the path plan.

{In contrast with kinodynamic motion planners and trajectory tracking approaches, we allow for Lagrangian systems to provably and safely fulfil path plans on the configuration space, rather than on the full state-space. We also reduce the offline computational burden of performing reachability analysis, and instead propose an abstraction of the obstacle-free space in the form of ellipsoids. Therefore, obstacles in the workspace are not explicitly inserted into the hybrid feedback control law proposed, retaining some convexity properties of the environment. Ellipsoids have already been used in the context of barrier functions in \cite{verginis2019closed} as a way of approximating the volume of an agent.}
\section{Preliminaries}
\label{sec:Preliminaries}
Real numbers are denoted by $\Real$, while $\Real^n$ is the $n$-dimensional real vector space.
A class $\mathcal{K}$ function is a function $\alpha(r)$, defined for $r \in [0,a)$, that is strictly increasing and with $\alpha(0) = 0$.
An extended class $\classKinf$ function is a class $\mathcal{K}$ function defined over $(-b,c)$, for some $b,c > 0$.
A function $\beta(r,s)$, defined for $r \in [0,a)$ and $s \in [0,\infty)$, belongs to class $\mathcal{KL}$ if $\beta(r,\cdot)$ belongs to class $\mathcal{K}$, and $\beta(\cdot, s)$ is decreasing, i.e. $\beta(\cdot, s) \to 0$ as $s \to \infty$.
$\intC$ and $\boundC$ denote the interior and boundary of the set $\setC$, respectively.
\fdsb{All proofs in this paper are provided in the appendix.}
An ellipsoid can be described by the following equation, where $A$ is a symmetric, positive definite matrix, $p \in \Real^{n}$ and $p_0 \in \Real^{n}$ its center:
\begin{equation}\label{eq:ellipse}
    (p - p_0)^T A (p - p_0) = 1.
\end{equation}

Let $x \in \Real^n$ and $u \in \Real^m$ be state and input of a nonlinear control-affine system
\begin{equation}
    \dot{x} = f(x) + g(x)u, \label{eq:system}
\end{equation}
with locally Lipschitz continuous functions $f : \Real^n \rightarrow \Real^n$ and $g : \Real^n \rightarrow \Real^{n\times m}$. 
%
%

The bounded subset $\mathcal{X} \subset \Real^{n}$ defines a workspace (configuration space), divided into a set of \emph{obstacles} $\mathcal{X}_{\text{obs}} \subset \mathcal{X}$ and the \emph{free space} $\mathcal{X}_{\text{free}} = \mathcal{X} \setminus \mathcal{X}_{\text{obs}}$.

\subsection{Path Plan}

\begin{defn}[Path Plan]\label{def:path}
    A path plan is a finite sequence of $N+1$ configurations $\mathbf{\bar{x}} = \mathbf{x}_0, \mathbf{x}_1, \dots, \mathbf{x}_N$, with $\mathbf{x}_0$ the initial configuration, $\mathbf{x}_i \in \mathcal{X}_{\text{free}}$ and $\mathbf{x}_N \in \mathcal{X}_{\text{goal}}$.
\end{defn}

\begin{assum}\label{assum:ellipse}
    Given a path plan $\mathbf{\bar{x}}$, there exist $N$ obstacle-free ellipsoids $\setC_i \subset \mathcal{X}_{\text{free}}$, {i.e. $N$ $A_i$ and $p_{0,i}$ as in \eqref{eq:ellipse},} such that $\mathbf{x}_i, \mathbf{x}_{i+1} \in \setC_i$, for $i = 0, \dots, N - 1$.
\end{assum}


    Assumption \ref{assum:ellipse} is guaranteed by the existence of a minimum clearance amongst the obstacles in the environment, such that a trajectory connecting two configurations is not unique. Minimum clearance is a realistic and natural assumption.

\begin{rem} A path plan can be obtained through different methods and algorithms. The choice of which to use is often problem specific.
By proposing the use of a generic definition of a path plan we allow the approach proposed in this paper to be applicable to a large variety of problems.
For instance, one can apply it to sampling-based motion planning algorithms, such as RRG, RRT$^\star$, PRM and PRM$^\star$ \cite{karaman2011sampling}, which build graphs (trees) on the configuration space of the system. 
\end{rem}

\subsection{Control Lyapunov Functions}


\begin{defn}[CLF]\label{def:esclf}
    \cite{khalil2002nonlinear} A continuously differentiable function $V: D \subset \Real^n \to \Real_{\ge 0}$ is an asymptotically stabilizing control Lyapunov function for \eqref{eq:system} if there exist continuous, positive definite functions $W_1(x), W_2(x), W_3(x)$ and {$u \in \Real^n$ such that the following hold for all $x \in D \subset \Real^n$}:
    \begin{align}
        & W_1(x) \le V(x) \le W_2(x) \\
        & \dot{V}(x) = \frac{\partial V(x)}{\partial x}f(x) + \frac{\partial V(x)}{\partial x}g(x)u \leq -W_3(x).
    \end{align}
\end{defn}

Definition \ref{def:esclf} allows one to define the set of stabilizing controllers for every $x \in D$ as
\begin{align*}
    \text{\small $K_{\text{CLF}}(x) := \left\{u \in \Real^m: \frac{\partial V(x)}{\partial x} \left(f(x) + g(x)u \right) \leq -W_3(x)\right\}.$}
\end{align*}
It follows that, if there exists a CLF for a nonlinear control-affine system \eqref{eq:system}, then \eqref{eq:system} is asymptotically stable under a control law $u(x)\in K_{\text{CLF}}(x)$ (under the assumption that $u(x)$ is locally Lipschitz continuous) \cite{ames2014rapidly}, and the solutions of the closed-loop system satisfy the following corollary.

\begin{cor} \label{def:beta_f}
    \cite[Thm. 4.9]{khalil2002nonlinear} Let $V(x)$ be a CLF for \eqref{eq:system} satisfying the conditions in Def.~\ref{def:esclf}. Moreover, for $r$ and $c$ chosen such that $B_r = \{\|x\| \le r\} \subset D$ and $c < \min_{\|x\| = r} W_1(x)$, then every trajectory {of \eqref{eq:system} under a control law $u(x)\in K_{\text{CLF}}(x)$} starting at $\{x \in B_r : W_2 \le c\}$ satisfies
    \begin{equation}
        \|x(t)\| \leq \beta(\|x(t_0)\|, t - t_0),
    \end{equation}
    for some class $\mathcal{KL}$ function $\beta$.
\end{cor}

This allows for the formulation of optimization-based controllers, without the need for explicitly defining a closed-form feedback control law $u(x)$ \cite{ames2019control}.

\subsection{Control Barrier Functions}
Consider a set $\setC$ defined as the superlevel set of a continuously differentiable function $h: \Real^n \to \Real$ as
\begin{align}
    \setC &= \left\{ x \in \Real^n: h(x) \ge 0 \right\} \label{eq:setCa}\\
    \boundC &= \left\{ x \in \Real^n: h(x) = 0 \right\}\label{eq:setCb}\\ 
    \intC &= \left\{ x \in \Real^n: h(x) > 0 \right\}.\label{eq:setCc}
\end{align}

\begin{defn}
    A set $\setC$ is forward invariant with respect to system \eqref{eq:system} if, for every initial condition $x(t_0) \in \setC$, its solutions remain within $\setC$, i.e. $x(t) \in \setC$ for $\forall t \geq t_0$. 
\end{defn}

The system \eqref{eq:system} is safe with respect to a set $\setC$ if such a set is forward invariant, allowing one to name $\setC$ a \emph{safe set}.

\begin{defn}[CBF]\label{def:cbf}
    \cite{ames2019control} Let $\setC \subseteq D \subset \Real^n$ be the superlevel set of a continuously differentiable function $h: \Real^n \to \Real$. Then $h$ is a control barrier function if there exists an extended class $\classKinf$ function $\alpha$ such that for the control system \eqref{eq:system}, the following holds $\forall x \in D$:
    \begin{equation}\label{eq:cbf_}
        \text{\small $\sup_{u \in \Real^m} \left[\frac{\partial h(x)}{\partial x} f(x) + \frac{\partial h(x)}{\partial x} g(x) u\right] \geq -\alpha(h(x)).$}
    \end{equation}
\end{defn}

Definition \ref{def:cbf} allows the definition of the set of control values that render $\setC$ forward invariant (i.e. safe):
\begin{align*}
    \text{\small $K_{\text{CBF}}(x) := \left\{u \in \Real^m: \frac{\partial h(x)}{\partial x} f(x) + \frac{\partial h(x)}{\partial x} g(x) u \geq -\alpha(h(x))\right\}.$}
\end{align*}
It follows that, if there exists a CBF for \eqref{eq:system}, then $\setC$ is forward invariant  with respect to system \eqref{eq:system} under a control law $u(x)\in K_{\text{CBF}}(x)$ (under the assumption that $u(x)$ is locally Lipschitz continuous) \cite[Cor. 2]{ames2017control}. 

\section{Problem Formulation}
\label{sec:Problem}
It is of utmost importance for safety-critical systems to consider the effects of integrating motion planning and closed-loop control and to provide formal guarantees on the trajectory of the system. The goal of this paper is to provide a way of connecting motion planning algorithms to feedback control laws in order to ensure safety and asymptotic stability of Lagrangian systems.

Lagrangian systems can be modelled as
\begin{equation}\label{eq:system2}
        \dot{p} = v \qquad \dot{v} = f(p,v) + g(p,v)u,
\end{equation}
where $p \in \Real^{n}$ is the \emph{configuration space}, $v \in \Real^{n}$ is the \emph{internal dynamics} space that governs the motion, and $f : \Real^{n} \times \Real^{n} \rightarrow \Real^n$ and {$g : \Real^n \times \Real^{n} \rightarrow \Real^{n\times n}$} are locally Lipschitz continuous functions. We require that $g(p,v)$ has full row rank for each $(p,v) \in \Real^{n} \times \Real^{n}$, hence modelling mechanical Lagrangian systems. Note that Lagrangian systems are modelled by a positive definite $g(p,v)$ \cite{vidyasagar1989robot} which is implied if $g(p,v)$ is square and full row rank. 

The problem is posed as follows.

\begin{prob} \label{prob:prob}
    Given a Lagrangian system \eqref{eq:system2}, a workspace $\mathcal{X} \subset \Real^{n}$ divided into obstacles $\mathcal{X}_{\text{obs}}$ and free space $\mathcal{X}_{\text{free}}$, and a path plan $\mathbf{\bar{x}}$, find a feedback control law $u(x)$ that drives the system from an initial configuration $\mathbf{x}_0 \in \mathcal{X}_{\text{free}}$ to a goal region $\mathcal{X}_{\text{goal}} \subset \mathcal{X}$ while ensuring safety (obstacle avoidance) and asymptotic stability.
\end{prob}

\section{Solution}
\label{sec:Solution}

The core idea to solve Prob.~\ref{prob:prob} is to explore the properties of both CLFs and CBFs in order to compose a sequential solution that guarantees safety and asymptotic stability. Given a path plan $\mathbf{\bar{x}} = \mathbf{x}_0, \mathbf{x}_1, \dots, \mathbf{x}_N$, we propose to find $N$ ellipsoids (convex sets) $\setC_i \subset \mathcal{X}_{\text{free}}$ enclosing each consecutive pair of $\mathbf{\bar{x}}$, i.e. each edge connecting $\mathbf{x}_{i}$ to $\mathbf{x}_{i+1}$, for $i = 0, \dots, N-1$. 
Several algorithms can be used in order to find obstacle-free ellipsoids in the workspace; to name one, we refer the reader to \cite{deits2015computing}.
The solution presented here renders these ellipsoids into \emph{safe sets} by constructing appropriate CBFs, designed to keep the system within its bounds, i.e. render $\setC_i$ controlled invariant.

Asymptotic stability is brought in by interpreting $\mathbf{\bar{x}}$ as a sequence of reachability problems, in which $N$ CLFs are written as to consecutively bring the system from the vicinity of $\mathbf{x}_i$ to the vicinity of $\mathbf{x}_{i+1}$, $i = 0, \dots, N-1$. In short, we build a sequence of reach-avoid control problems, where the dynamical system is supposed to remain within $\setC_{i+1}$ while being controlled from $\mathbf{x}_i$ to $\mathbf{x}_{i+1}$.

Although convergence to each $\mathbf{x}_{i}$ happens asymptotically, finite time convergence is guaranteed to the vicinity of $\mathbf{x}_{i}$ due to the use of CLFs. 
This motivates us to propose a switching mechanism between QPs: since $\mathbf{x}_{i}$ has to be inside both $\setC_{i-1}$ and $\setC_{i}$, for $i = 1, \dots, N-1$, and if the intersection of such sets is not trivial, one can switch in a provably-correct manner to the next control law (i.e. reach $\mathbf{x}_{i+1}$) whenever the system reaches $\setC_i \cap \setC_{i+1}$ with a proper velocity, which happens in finite time.

\subsection{Deriving a CBF from an ellipsoid}
Recall from Def.~\ref{def:path} and Asm.~\ref{assum:ellipse} that, given a path plan $\mathbf{\bar{x}}$, there exist obstacle-free ellipsoids $\setC_i \subset \mathcal{X}_{\text{free}}$ enclosing $\mathbf{x}_i, \mathbf{x}_{i+1}$.
The objective of this section is to ensure each $\setC_i$ is a \emph{valid safe set}  for Lagrangian systems \eqref{eq:system2}. 
We aim at defining the set of control values that render each $\setC_i$ forward invariant to~\eqref{eq:system2}.
For simplicity, we drop the index $i$ in the rest of this subsection.
Let us start by defining $\setC$ on the workspace of Prob.~\ref{prob:prob}:
\begin{equation}\label{eq:setC}
    \setC := \left\{ p \in \mathcal{X}_{\text{free}} \subset \Real^{n}: h(p) \ge 0 \right\}.
\end{equation}

Since $\setC$ is a superlevel set of $h: \Real^{n} \to \Real$, and $\setC$ is supposed to encode obstacle-free ellipsoids in the workspace, we propose to define $h(p)$ directly from the generic equation of ellipsoids, i.e. as
\begin{equation}\label{eq:h}
    h(p) := 1 - (p - p_0)^T A (p - p_0),
\end{equation}
where $p_0 \in \mathcal{X}_{\text{free}} \subset \Real^{n}$ is the center of the ellipsoid. Note that \eqref{eq:h} defines $\setC$ in accordance with \eqref{eq:setCa}-\eqref{eq:setCc}.
Such definition does not necessarily satisfy Definition~\ref{def:cbf}, i.e. \eqref{eq:cbf_} may not hold,  because the system in \eqref{eq:system2} is a relative degree two system. In other words, $\frac{\partial h(p)}{\partial p} \dot{p} = \frac{\partial h(p)}{\partial p} v$ so that no control input $u$ appears in \eqref{eq:cbf_}, which evaluates in this case to $\frac{\partial h(p)}{\partial p} v\ge -\alpha(h(p))$, where $\alpha$ is a user-defined extended class $\classKinf$ function. A solution, known in the literature as higher-order CBF \cite{nguyen2016exponential,xiao2019control}, is to introduce new CBFs that depend on partial derivatives of $h(p)$. For a relative degree two system, only one CBF needs to be introduced, here denoted by $h': \Real^{n} \times \Real^{n} \to \Real$ and defined as 
\begin{align}\label{eq:hprime}
    h'(p,v) &:= \frac{\partial h(p)}{\partial p} v + \alpha(h(p)) \\
            &= - 2 (p - p_0)^T A v + \alpha(h(p)).
\end{align}
 Naturally, we associate
\begin{equation}\label{eq:setCprime}
    \setC' := \left\{ (p,v) \in \Real^{n} \times \Real^{n}: h'(p,v) \ge 0 \right\},
\end{equation}
with $h'(p,v)$ that now depends on both $p$ and $v$. For an extended class $\classKinf$ function $\alpha'$ and if
\begin{equation}\label{eq:Kcbf}
        \text{\small $\frac{\partial h'(p,v)}{\partial p} v + \frac{\partial h'(p,v)}{\partial v} \left(f(p,v) + g(p,v)u\right) \geq  - \alpha'(h'(p,v))$}
\end{equation}
for all $p \in \setD\supseteq \setC$ and $(p,v) \in \setD'\supseteq\setC'$ for some open sets $\setD$ and $\setD'$, it follows that 
$h'(p(t),v(t))\ge 0$ if $h'(p(0),v(0))\ge 0$. This implies that $h(p(t))\ge 0$ if $h(p(0))\ge 0$. Note that \eqref{eq:Kcbf} results from applying Definition~\ref{def:cbf} to $h'(p,v)$.
%

\begin{lem}\label{th:cbf}
    Let $h(p)$ and $h'(p,v)$ be defined as in \eqref{eq:h} and \eqref{eq:hprime}, respectively, and let $\setC$ and $\setC'$ be the corresponding sets as defined in \eqref{eq:setC} and \eqref{eq:setCprime}. 
    Then $h'(p,v)$ is a feasible control barrier function for almost all $(p,v)$ for Lagrangian systems \eqref{eq:system2}, i.e. there exist $u(p,v)$ and extended class $\classKinf$ functions $\alpha, \alpha'$ such that
    \begin{align}\label{eq:KcbfEl}
        \small
        &v^T \frac{\partial^2 h(p)}{\partial p^2} v + \frac{\partial \alpha(h(p))}{\partial p} v + \frac{\partial h(p)}{\partial p} f(p,v) \nonumber\\
        & \hspace{2cm} + \frac{\partial h(p)}{\partial p} g(p,v)u + \alpha'(h'(p,v)) \geq 0
    \end{align}
    holds for all $p \in \setD \supseteq \setC$ and $(p,v) \in \setD' \supseteq \setC'$, except for $p = p_0$, where it is infeasible if $v$ is unbounded. {In the case of $v$ being bounded, then $h'(p,v)$ is also feasible at $p = p_0$.}
\end{lem}


\begin{rem}\label{rem:remark}
    The CLF proposed in Sec.~\ref{sec:CLF} enforces the states of the system to be bounded, therefore making the CBF proposed in this section a valid one.
\end{rem}

Following these results, the set of control values that render $\setC'$ forward invariant, and therefore $\setC$ as well, is defined as follows:
\begin{align}\label{eq:K_cbf}
    &K_{\text{CBF}}(p,v) := \left\{u \in \Real^m: -2v^TAv + \frac{\partial \alpha(h(p))}{\partial p} v \right.  \\
    &\left.- 2(p-p_0)^TA (f(p,v) + g(p,v)u) + \alpha'(h'(p,v)) \ge 0\right\}.\nonumber
\end{align}

\begin{figure}
    \centering
    \includegraphics[width=0.9\linewidth]{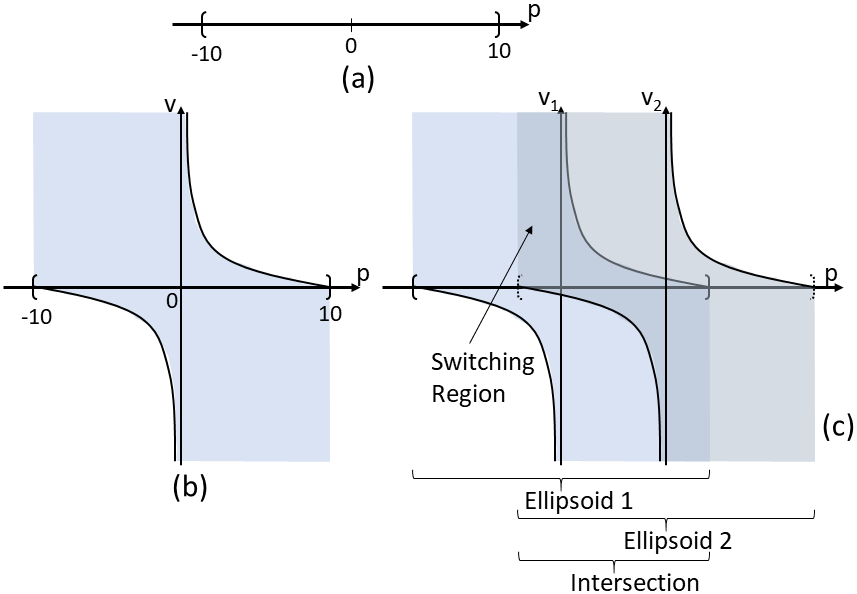}
    \vspace{-3mm}
    \caption{Illustration for Ex.~\ref{ex:one} and \ref{ex:two}. (a) a one-dimensional ellipsoid of radius 10, which corresponds to $\setC$. (b) depicts, in light blue, $\setC'$ for an arbitrary $\alpha(h(p))$. (c) two one-dimensional ellipsoids $\setC_1$ and $\setC_2$ are depicted, whose centers are aligned with the axes $v_1$ and $v_2$, respectively; regions in light blue and gray correspond to $\setC'_1$ and $\setC'_2$, respectively; the overlapping region, pointed in the figure as Switching Region, is the subset of $\Real^2$ such that for every $p$ and $v$ inside it, $p \in \setC_1 \cap \setC_2$, and $(p,v) \in \setC'_1 \cap \setC'_2$.}
    \label{fig:velocity}
\end{figure}

\begin{exmp}\label{ex:one}
    For clarification, let us study a one-dimensional example. Suppose there exists an ellipsoid centered at $p_0 = 0$ with $A = 0.01$, i.e. its radius is 10, as depicted in Fig.~\ref{fig:velocity}a. The safe-set $\setC$ rendered by such ellipsoid is as follows $$\setC = \left\{ p \in \Real: 1 - A p^2 \ge 0 \right\},$$ i.e. $\setC = \left\{ p \in [-10,10]\right\}$. Similarly, we define $$\setC' = \left\{ (p,v) \in \Real^2: -2pAv + \alpha(h(p)) \ge 0\right\}$$ which can be used to find admissible values of $v$ for each point $p$ and with respect to feasibility of \eqref{eq:KcbfEl}. Fig.~\ref{fig:velocity}b shows $\setC'$ for an arbitrary $\alpha(h(p))$.
\end{exmp}

\subsection{Valid CLFs and QP}\label{sec:CLF}
The purpose of this subsection is to define the structure of the $N$ feedback controllers that ensure asymptotic stability and safety of the system throughout the path plan $\mathbf{\bar{x}}$. We start by defining the CLF to each $\mathbf{x_i} \in \mathbf{\bar{x}}$, and then its formulation into a QP alongside the corresponding CBF.

Let us define $\boldsymbol{p}_i$, with $i = 1, \dots, N$, as the error between the position of the system \eqref{eq:system2} and the current goal configuration $\mathbf{x}_i$, i.e. $\boldsymbol{p}_i := \mathbf{x}_i - p$. We propose a candidate CLF $V(\boldsymbol{p}_i,v)$ to the Lagrangian system \eqref{eq:system2} such that its equilibrium state is where the position error and velocity are both zero, i.e. $V(\boldsymbol{p}_i,v) = 0$ if and only if $\boldsymbol{p}_i = 0$ and $v = 0$. More specifically, we propose the following candidate function
\begin{align} \label{eq:clf_P}
    V(\boldsymbol{p}_i,v) &:= \frac{1}{2}\begin{bmatrix}\boldsymbol{p}_i^T & v^T \end{bmatrix} P \begin{bmatrix}\boldsymbol{p}_i^T & v^T \end{bmatrix}^T \\
    &= \frac{1}{2}\begin{bmatrix}\boldsymbol{p}_i^T & v^T \end{bmatrix} \begin{bmatrix}P_1 & P_2 \\ P_2^T & P_3\end{bmatrix} \begin{bmatrix}\boldsymbol{p}_i^T & v^T \end{bmatrix}^T
\end{align}
with $P$ a symmetric, positive definite block matrix, i.e. $P_1$ and $P_3 - P_2^T P_1^{-1} P_2$ are positive definite, and $P_2$ negative definite (follows due to properties of the Schur complement). Such a CLF candidate is a common choice for Lagrangian systems, as it can represent the system's energy.

\begin{lem}\label{lem:clf_valid}
    Consider the Lagrangian system \eqref{eq:system2}. The candidate Lyapunov function $V(\boldsymbol{p}_i,v)$ in \eqref{eq:clf_P} is a valid CLF, and therefore renders the origin of the system (globally) asymptotically stable.
\end{lem}

%
    %
%
%

Lemma \ref{lem:clf_valid} allows us to define the set of asymptotic stabilizing control values as follows:
\begin{align}\label{eq:K_clf}
    K_{\text{CLF}}(p,v) &:= \left\{u \in \Real^m: -\boldsymbol{p}_i^T P_1 v - v^T P_2^T v \right.  \\
    &\hspace{-1.7cm}\left. + (\boldsymbol{p}_i^T P_2 + v^T P_3) (f(p,v) + g(p,v)u) \le - W_3(\boldsymbol{p}_i,v)\right\} \nonumber.
\end{align}

Following what has been proposed in recent related works (e.g. \cite{ames2017control,ames2019control}), we also consider a feedback control law in the form of a quadratic optimization problem. From a path plan $\mathbf{\bar{x}}$ of size $N + 1$, and having the Lagrangian system starting at $p(t_0) = \mathbf{x}_0$ and $v(t_0) = 0$, we write a QP that drives the system asymptotically to $p(t) = \mathbf{x}_i$ and $v(t) = 0$, $t > t_0$, while ensuring it stays within the corresponding safe set $\setC_i$, for $i = 1, \dots, N$. The control law is the following
\begin{align}\label{eq:qp1}
    u (p,v) &= \argmin_{u \in \Real^{m}} \quad \frac{1}{2} u^T H(p,v)u \\
    \text{s.t. }\; & u \in K_{\text{CLF}}(p,v)\\
    & u \in K_{\text{CBF}}(p,v),
\end{align}
%
%
where $H(p,v)$ is a user-defined positive definite matrix. 

\begin{assum}
    The formulation of the QP \eqref{eq:qp1} generates Lipschitz continuous controllers.
\end{assum}

{Sufficient conditions under which such a control law is Lipschitz continuous are discussed in \cite{morris2015continuity}. Two more ways to prove Lipschitz continuity of combined CLF/CBF control laws are presented in the proofs of \cite[Theorem 11]{xu2015robustness} (based on the KKT conditions) and in \cite[Theorem 1]{jankovic2018robust}.}

\begin{thm}\label{th:main}
    Consider a system~\eqref{eq:system2} with initial states $p_0,v_0$ inside the safe sets $\setC$ and $\setC'$ defined by CBFs \eqref{eq:h} and \eqref{eq:hprime}, i.e., $p_0\in\setC$ and  $(p_0,v_0)\in\setC'$. Also consider $V(\boldsymbol{p}_i,v)$ in \eqref{eq:clf_P}, with $\boldsymbol{p}_i = \mathbf{x}_i - p$ and $\mathbf{x}_i \in \intC$. Then the QP \eqref{eq:qp1} is always feasible by a suitable choice of $\alpha$, $\alpha'$ and $W_3$. 
\end{thm}


{In short, we prove in Theorem \ref{th:cbf} that the QP proposed in \eqref{eq:qp1} is always feasible, i.e. there always exists a control input $u \in \Real^n$ that satisfies the constraints imposed by both CBF and CLF. Therefore, the trajectories of the system will always remain within the safe sets $\setC$ and $\setC'$ while asymptotically converging to the next location according to CLF $V$.}

\subsection{Switching Mechanism}
In the previous sections we defined our approach to designing CBFs from ellipsoids, which ensure the safety of the system, and CLFs, which ensure asymptotic convergence to each configuration of a path plan. We also defined $N$ QPs that connect the $i$-th CBF and CLF into one feedback control problem. Now, the last step of the proposed solution is to describe a correct switching mechanism so that the controller is well-defined even when transitioning from the $i$-th to the $(i+1)$-th QP. Note that in the sense of the proposed CLF formulation, switching from $V(\boldsymbol{p}_i,v)$ to $V(\boldsymbol{p}_{i+1},v)$ is a simple coordinate transformation on the equilibrium point, and therefore stability is not lost.

Regarding safety, it is desired that the system remains safe during all time, including when it switches from one safe set to another. In the scope of our paper, this means that in order for the system to switch from the $i$-th QP to the next one, it must hold that $p \in \setC_i \cap \setC_{i+1}$ and $(p,v) \in \setC'_i \cap \setC'_{i+1}$. Since we motivate our solution by claiming that it is enough to build the CBF directly from geometrical constraints, we must now prove that if {$\setC_i$ and $\setC_{i+1}$ intersect, i.e. $\setC_i \cap \setC_{i+1} \neq \emptyset$, then so do $\setC'_i$ and $\setC'_{i+1}$, i.e. $\setC'_i \cap \setC'_{i+1} \neq \emptyset$.}

\begin{thm}\label{th:switch}
    Let $h(p)$ as defined in \eqref{eq:h} be a valid CBF for the Lagrangian system \eqref{eq:system2}. Also let $\setC$ and $\setC'$ be defined as in \eqref{eq:setC} and \eqref{eq:setCprime}, respectively. For two different sets $\setC_i$ and $\setC_j$, it holds that if $\setC_i \cap \setC_j \neq \emptyset$, then $\setC'_i \cap \setC'_j \neq \emptyset$.
\end{thm}


\begin{rem}
    Note that we design our solution such that every $\mathbf{x}_i \in \mathbf{\bar{x}}$, for $i = 1, \dots, N-1$, is inside $C_{i-1} \cap C_{i}$, and that we design CLFs to asymptotically bring the system to $\mathbf{x}_i$, therefore to the switching region $C'_{i-1} \cap C'_{i}$.
\end{rem}

\begin{exmp}\label{ex:two}
    Continuing the previous one-dimensional example, let us suppose we have two ellipsoids, $\setC_1$ and $\setC_2$, that intersect each other, as depicted in Fig.~\ref{fig:velocity}(c). The regions in blue and gray are, respectively, $\setC'_1$ and $\setC'_2$.
    %
\end{exmp}

\subsection{Discussion}
Even though the proposed CBF generates a safe set $\setC'$ that is unbounded on $v$ (see \eqref{eq:hprime}, especially for $p = p_0$), our approach also uses a CLF which, besides being asymptotic stabilizing, generates invariant sublevel sets on the state-space. Therefore, the bounds on $v$ are determined by the state of the system at each switch between safe sets.

Although the proposed CLF guarantees asymptotic convergence to each $\mathbf{x}_i \in \mathbf{\bar{x}}$ as stated in Def.~\ref{def:esclf}, we prove in Theorem~\ref{th:switch} that the intersection of two consecutive safe sets $\setC'_i, \setC'_{i+1}$ is not a singleton. Therefore, there exists an $\epsilon$-ball $B_{\epsilon}$ centered on $(\mathbf{x}_i, v=0)$ such that $B_{\epsilon} \subset \setC'_i \cap \setC'_{i+1}$. As a consequence of Theorem~\ref{th:cbf}, there exists a class $\mathcal{KL}$ function $\beta$ as in Cor.~\ref{def:beta_f} that describes the trajectory of the system, which intersects $B_{\epsilon}$ in finite time.

Lastly, Theorems \ref{th:cbf} and \ref{th:switch} together render the proposed approach sound and complete. In other words, the system will always perform a trajectory that is guaranteed to be safe, i.e. to remain within the sequence of $\setC, \setC'$, and that (asymptotically) reaches a configuration in the goal region. 

\section{Case Studies} \label{sec:cases}
\fdsb{We demonstrate the application of the proposed approach in two case studies, both implemented in MATLAB. We use $\alpha(h(p)) = k_1 h^3(p)$ and $\alpha'(h'(p,v)) = k_2 h'^3(p,v)$ in both case studies, with $k_1,k_2 >0$ tunable parameters. The sequences of ellipsoids and the path plans in both scenarios were manually designed and placed by us. Implementation of a tool for automatically planning such a sequence of obstacle-free ellipsoids is out of the scope of this paper and left as an interesting direction for future work.}

\begin{figure}
    \centering
    \includegraphics[trim={29.2mm 54mm 18mm 49mm},clip, width=\linewidth]{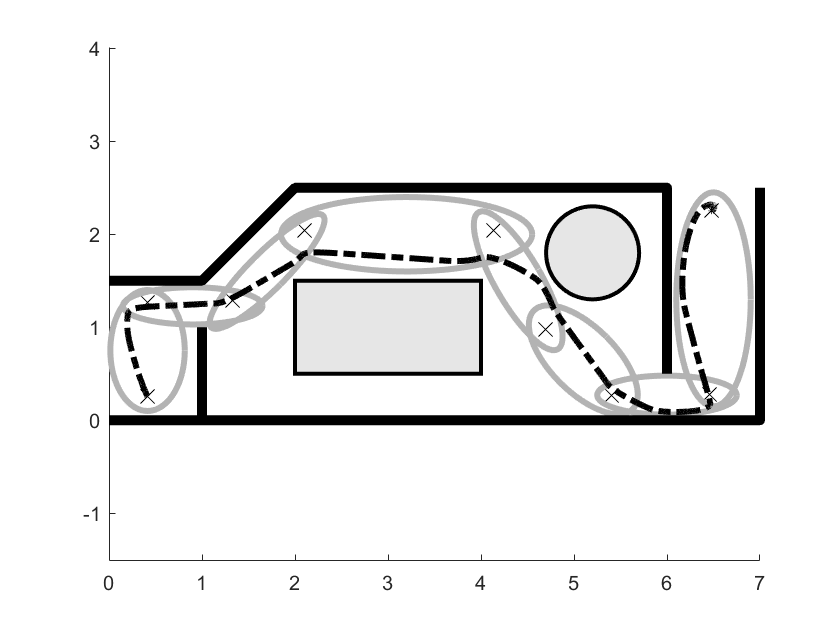}
    \vspace{-5mm}
    \caption{\fdsb{Simulation of our proposed in a 2D office-like environment. Black thick lines represent walls, and in gray are two obstacles. The path plan is in black `x' markers, and the obstacle-free ellipsoids in gray.}}
    \label{fig:DoubleIntegratorEx}
\end{figure}

\fdsb{The first case study is in a two-dimensional environment, $p~:=~ \begin{bmatrix}p_x & p_y \end{bmatrix}^T \in \Real^2$ and $v := \begin{bmatrix}v_x & v_y \end{bmatrix}^T \in \Real^2$, with a system with damped and coupled dynamics as $\dot{p}~=~v$, $m\dot{v}_x = -b v_x + u_x - 0.5u_y$ and $m\dot{v}_y = -b v_y + u_y$, where $m,b > 0$ are the mass and damping parameters, respectively. Fig.~\ref{fig:DoubleIntegratorEx} shows that our approach successfully prevents the system from leaving the obstacle-free ellipsoids while also ensuring it converges to the goal position.}
%
%



\begin{figure*}
    \centering
    \subfloat[View from initial configuration\label{3d3a}]{\includegraphics[width=0.38\linewidth]{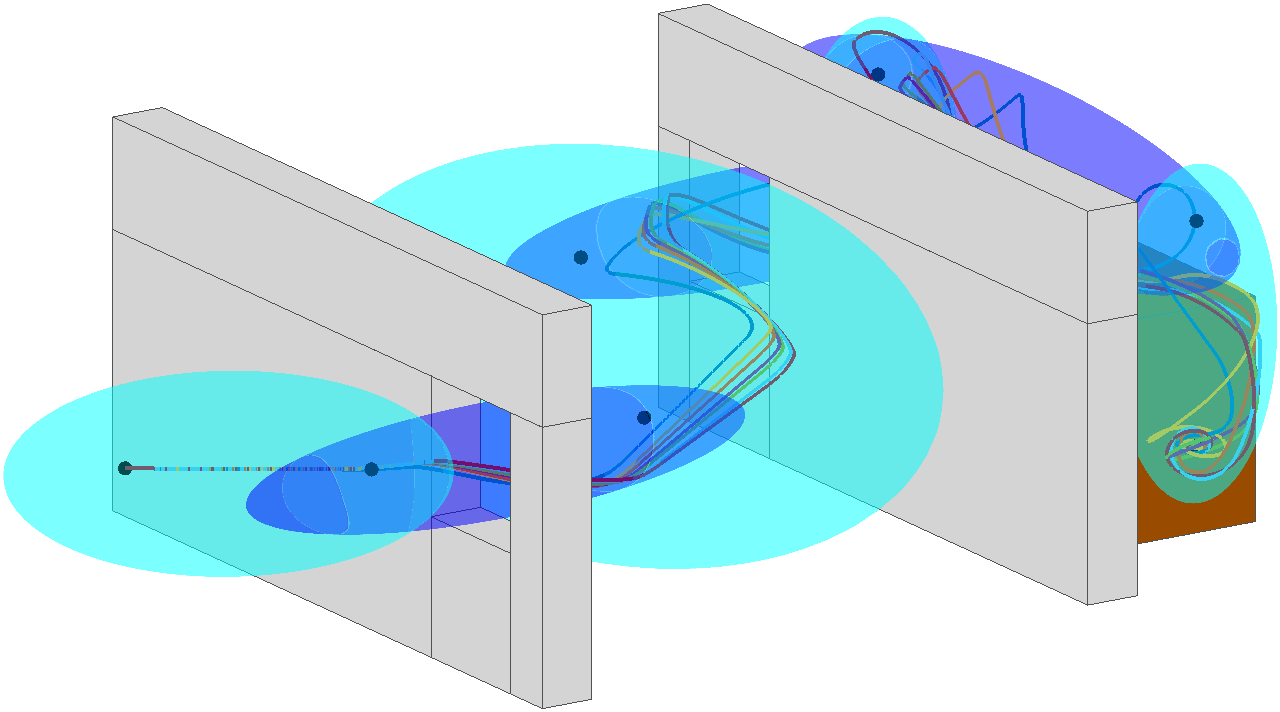}}
    \subfloat[View from goal configuration\label{3d3b}]{\includegraphics[width=0.3\linewidth]{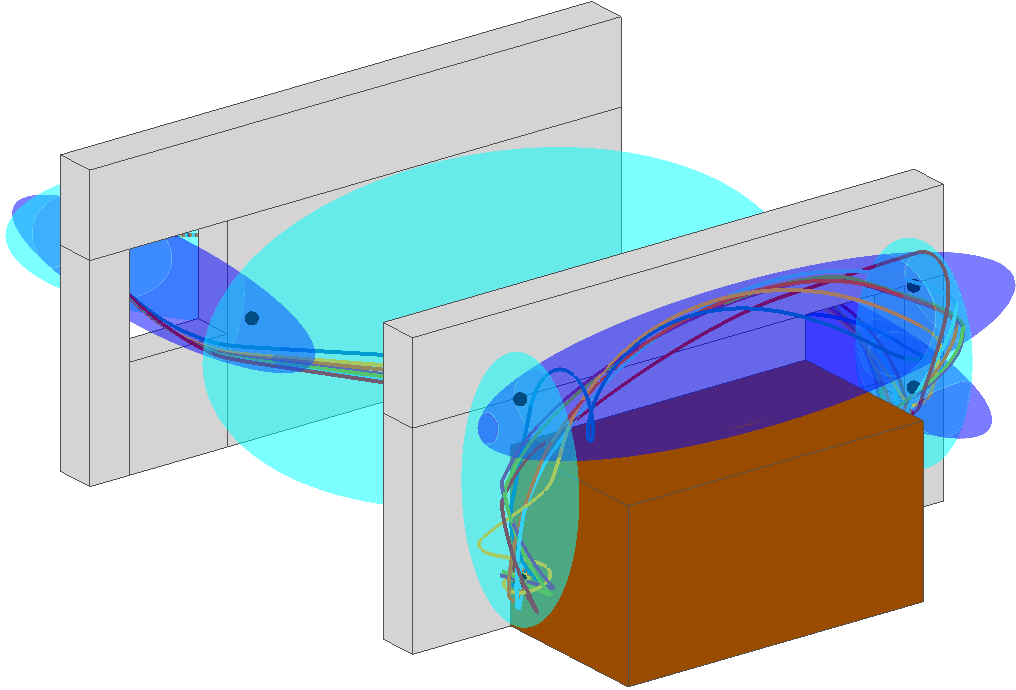}}
    \subfloat[Top-view\label{3d3c}]{\includegraphics[width=0.32\linewidth]{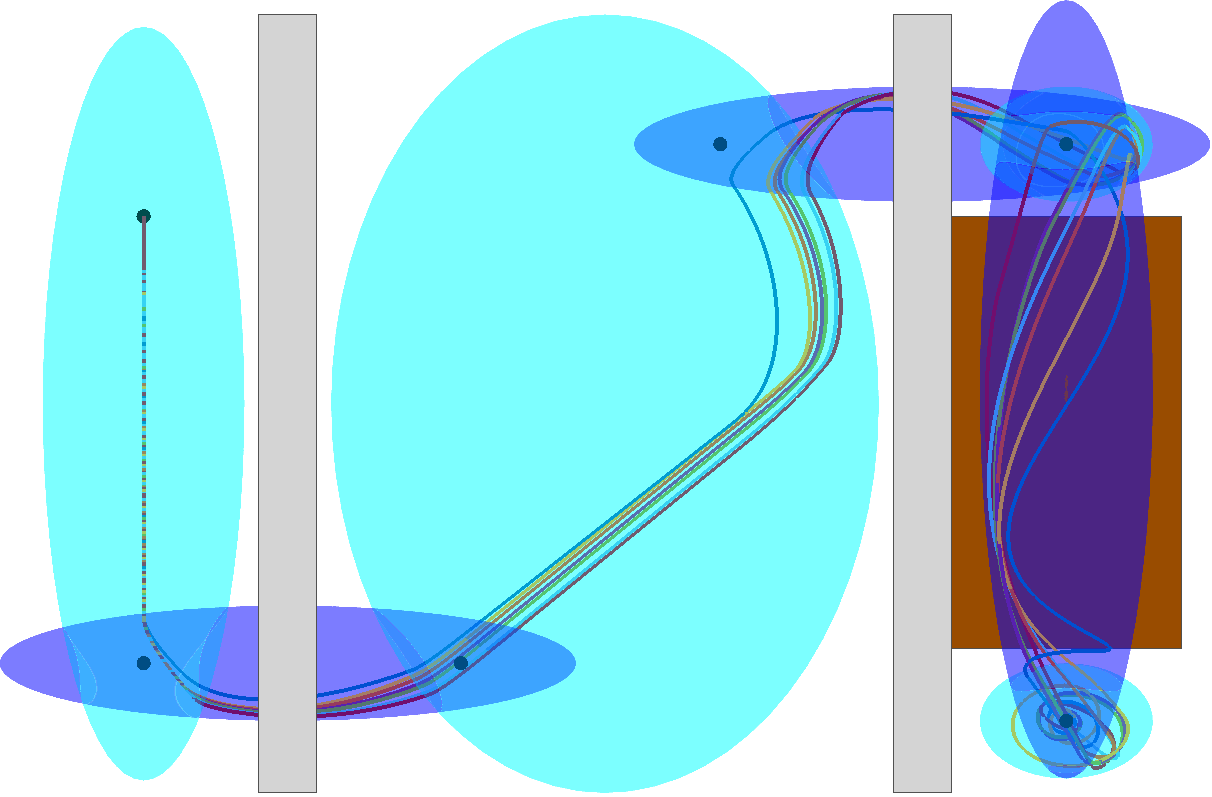}}
    \caption{Application of our proposed approach in a complex 3D environment, with two walls in gray and a box in brown. The path plan is displayed in black dots and the obstacle-free ellipsoids in light and dark blue. Several trajectories are displayed, in different colors, to highlight how the parameters in the CBF can influence the final result.}
    \label{fig:3D}
\end{figure*}

\fdsb{In the second case study we use a 3D single-integrator system, i.e. $\dot{p} = v$ and $\dot{v} = u$, with $(p,v) \in \Real^3\times\Real^3$. In this environment, shown in Fig.~\ref{fig:3D}, the system is required to go through holes in two walls and then over a box, before lowering its altitude. In order to show how the parameters of the CBF influence the resulting trajectory of the system, we run several simulations with different values of $k_1, k_2$ and display the collection of trajectories in Fig.~\ref{fig:3D}. Note how our approach guarantees the system stays within the safe sets, even when switching from large to relatively small ellipsoids, e.g. for going through the holes in the walls.}

The case studies highlight \fdsb{the modularity} of our approach and how it can be seen as an alternative to kinodynamic planners. Instead of finding an optimal trajectory and employing a trajectory-tracking controller, our approach explores as much of the obstacle-free space as possible, not constraining the system to tracking a prescribed trajectory \fdsb{and simply requiring a path plan and safe ellipsoidal regions. Differently from other approaches to collision avoidance, our approach does not require any special treatment of obstacles, such as abstractions or explicit inclusion in the control law. Instead, we choose to abstract the empty-space.} 

\section{Conclusion}
\label{sec:Conclusion}
We derived a provably safe, hybrid, closed-loop navigation algorithm for Lagrangian systems. We proposed feasible CBFs and CLFs, together with a switching mechanism, that allow such systems to safely follow path plans created on the configuration space. Throughout the paper the workspace is assumed to be known a priori. However, this assumption does not limit the applicability of the approach, which can be easily adapted to receding-horizon algorithms in environments mapped in real time.

Future work includes the consideration of dynamic obstacles in the workspace, \fdsb{inclusion of velocity and input bounds}, as well as expansion to multi-agent systems, such as formation navigation. Furthermore, we also plan to include temporal logics as a way of specifying desired motion preferences.

\section*{APPENDIX}
\emph{Proof of Lemma \ref{th:cbf}:}
Recall \eqref{eq:Kcbf} and take the derivative of \eqref{eq:hprime} with respect to $p$. Knowing that $\frac{\partial h'(p,v)}{\partial v} = \frac{\partial h(p)}{\partial p}$, we reach \eqref{eq:KcbfEl}. Note that \eqref{eq:KcbfEl} can always be made feasible by selecting an appropriate $u$ if $\frac{\partial h(p)}{\partial p}g(p,v) \neq 0$. In case such condition does not hold, note that $\frac{\partial h(p)}{\partial p}g(p,v) = 0$ if and only if $\frac{\partial h(p)}{\partial p} = 0$, since $g(p,v)$ has full row rank. Note further that $\frac{\partial h(p)}{\partial p} = 0$ if and only if $p = p_0$ since $h(p)$ is concave.
In order to ensure satisfaction of \eqref{eq:KcbfEl} in this singularity point ($p = p_0$), the following equation has to hold:
\begin{align}\label{eq:KcbfEl_sing}
    \text{\small $v^T \frac{\partial^2 h(p_0)}{\partial p^2} v + \frac{\partial \alpha(h(p_0))}{\partial p} v + \alpha'(h'(p_0,v)) \geq 0.$}
\end{align}
Analyzing this equation term by term, we know that (i) $\frac{\partial^2 h(p_0)}{\partial p^2} = -2A$ is negative definite, thus $v^T \frac{\partial^2 h(p_0)}{\partial p^2} v < 0$; (ii) $h'(p_0,v) = \alpha(h(p_0)) = \alpha(1) > 0$, and therefore $\alpha'(h'(p_0,v)) = \alpha'(\alpha(1)) > 0$; and (iii) the middle term is
\begin{align}
    \text{\small $\frac{\partial \alpha(h(p_0))}{\partial p}$} & \text{\small $= \frac{d \alpha(h(p_0))}{d h(p)} \frac{d h(p_0)}{d p}$} \label{eq:alpha_der}\\
    & \text{\small $= - 2 \frac{d \alpha(h(p_0))}{d h(p)} (p_0 - p_0)^T A = 0.$}
\end{align}
Substituting the previous terms into \eqref{eq:KcbfEl_sing} results in $-2v^TAv + \alpha'(\alpha(1)) \ge 0$, which cannot be satisfied if $v$ is unbounded. Therefore, $h'(p,v)$ is a valid CBF for all $(p,v) \in \setD' \supseteq \setC'$ except for $p = p_0$, where it is infeasible if $v$ is unbounded. {It is straightforward to see that if $v$ is bounded, then there exist class $\classKinf$ functions $\alpha, \alpha'$ such that $\alpha'(\alpha(1)) \ge 2v^TAv$.}
%
\qed

\emph{Proof of Lemma \ref{lem:clf_valid}:}
Note that $V(\boldsymbol{p}_i,v)$ is a positive definite function, and that according to Def.~\ref{def:esclf}, $\dot{V}(\boldsymbol{p}_i,v)$ must be smaller than or equal to a negative definite function $-W_3(\boldsymbol{p}_i,v)$. Taking the derivative of \eqref{eq:clf_P}
\begin{align}
    \dot{V} &= \begin{bmatrix}\boldsymbol{p}_i^T & v^T \end{bmatrix} P \begin{bmatrix}-v \\ f(p,v) + g(p,v)u \end{bmatrix} \\
    &= -\boldsymbol{p}_i^T P_1 v - v^T P_2^T v + (\boldsymbol{p}_i^T P_2 + v^T P_3)f(p,v) \nonumber \\
    & \hspace{1cm} + (\boldsymbol{p}_i^T P_2 + v^T P_3) g(p,v)u \le - W_3(\boldsymbol{p}_i,v) \label{eq:lyap_proof},
\end{align}
it holds that {there will always exist $u$ that satisfies} \eqref{eq:lyap_proof} for the case $(\boldsymbol{p}_i^T P_2 + v^T P_3) g(p,v) \neq 0$. When this condition does not hold, note that since $g(p,v)$ has full row rank, $(\boldsymbol{p}_i^T P_2 + v^T P_3) g(p,v) = 0$ if and only if $\boldsymbol{p}_i^T = -v^T P_3 P_2^{-1}$. Substituting this in \eqref{eq:lyap_proof} yields 
\begin{equation}
    v^T (P_3 P_2^{-1} P_1 - P_2^T) v + W_3(\boldsymbol{p}_i,v) \le 0, \label{eq:lyap_proof2}
\end{equation}
for which it is always possible to select a positive definite function $W_3(\boldsymbol{p}_i,v)$ that renders \eqref{eq:lyap_proof2} feasible if $P_3 P_2^{-1} P_1~-~P_2^T$ is negative definite.
Recall that we define $P$ as a symmetric, positive definite, block matrix, implying that $P_1$ and $P_3 - P_2^T P_1^{-1} P_2$ are both positive definite; we also define $P_2$ negative definite.
Note that all these constraints in the matrices forming $P$ are not conflicting. For instance, it is always possible to choose $P_1,P_2,P_3$ diagonal, such that $(P_3 - P_2^T P_1^{-1} P_2) P_2^{-1}P_1 = P_3 P_2^{-1} P_1 - P_2^T$ is negative definite.
Finally, see that $v=0$ could render \eqref{eq:lyap_proof2} infeasible. However, recall that the condition to enter \eqref{eq:lyap_proof2} is to have $\boldsymbol{p}_i^T = -v^T P_3 P_2^{-1}$. Therefore, if $v=0$, then $\boldsymbol{p}_i^T=0$, which implies $W_3(0,0) = 0$.
With this we finish the proof.
\qed

\emph{Proof of Theorem \ref{th:main}:}
Since $g(p,v)$ is positive definite and $f(p,v)$ is known, one can design a feasible (not necessarily optimal) controller $u = g(p,v)^T(g(p,v)g(p,v)^T)^{-1} (-f(p,v) + u^*)$ such that \eqref{eq:system2} becomes a double integrator.  The $K_{\text{CBF}}, K_{\text{CLF}}$ conditions in \eqref{eq:K_cbf},\eqref{eq:K_clf} can be respectively reorganized as follows:
\begin{align}
    & \text{\small $2(p - p_0)^T Au^* \le -2v^{T}Av + \frac{\partial \alpha(h(p))}{\partial p} v + \alpha'(h'(p,v))$} \label{eq:CBF_u} \\ 
    & \text{\small $(\boldsymbol{p}_i^T P_2 + v^T P_3) u^* \le v^T P_2^T v + \boldsymbol{p}_i^T P_1 v  - W_3(\boldsymbol{p}_i,v)$}\label{eq:CLF_u}.
\end{align}
where $w_1(p) := 2(p - p_0)^T A$ and $w_2(p,v) := \boldsymbol{p}_i^T P_2 + v^T P_3$. Let $c\ge 0$ be such that $V(\boldsymbol{p}_i(0),v_0)=c$ where $\boldsymbol{p}_i(0)=\mathbf{x}_i - p_0$. Define the corresponding sublevel set as $\mathcal{V}_i:=\{(p,v)\in\mathbb{R}^{n}\times\mathbb{R}^{n}|V(\mathbf{x}_i - p,v)\le c\}$ and note that $\mathcal{V}_i$ is compact. Since $v$ is bounded on this set, we can select $\alpha'$ such that \eqref{eq:CBF_u} by itself, i.e., when not considering \eqref{eq:CLF_u}, is always feasible on $\mathcal{V}_i$, as remarked in Rem.~\ref{rem:remark}. Note now that the case $w_1(p)=0$ and $w_2(p,v)\neq 0$ results in \eqref{eq:qp1} being feasible since the CBF condition holds by default as analyzed in Lemma \ref{th:cbf}. Subsequently, $u^*$ can be selected purely based on the CLF condition. For the case where $w_1(p)\neq 0$ and $w_2(p,v)= 0$,  the CLF condition holds by default as analyzed in Lemma \ref{lem:clf_valid}. Then $u^*$ can be selected purely based on the CBF condition. The case $w_1(p)=0$ and $w_2(p,v)= 0$ holds trivially.

Let us in the remainder focus on the combined case where $w_1(p)\neq 0$ and $w_2(p,v)\neq 0$. Note that the feasible region dictated \eqref{eq:CBF_u} and \eqref{eq:CLF_u} is given by the intersection of  two closed half-spaces with $w_1(p)$ and $w_2(p,v)$ the corresponding normal vectors perpendicular to the supporting hyperplanes, respectively. For $(p,v)$ such that $w_1(p)\times w_2(p,v)\neq 0$, i.e., $w_1(p)$ and $w_2(p,v)$ are not parallel, \eqref{eq:CBF_u} and \eqref{eq:CLF_u} are feasible simultaneously, even if $(p,v)$ are also in $\boundC'$, where $\alpha'(h'(p,v)) = 0$. On the other hand, if $w_1(p)$ and $w_2(p,v)$ are parallel, we have two cases to analyze: one in which they point towards the same direction, and one in which they point to opposite directions. It is straightforward to see that the former case makes \eqref{eq:CBF_u} and \eqref{eq:CLF_u} always feasible simultaneously. For the latter case, first note that due to the definition of our solution, every $p \in \boundC$, where $\alpha(h(p)) = 0$ and $h(p)=0$, is excluded since $w_1(p), w_2(p,v)$ cannot be parallel and point to opposite directions at the boundary of the safe set. Similar analysis holds for a neighborhood around $\boundC$. We have $\frac{w_1(p)}{\|w_1(p)\|}=-\frac{w_2(p,v)}{\|w_2(p,v)\|}$, with $\|w_1(p)\|$ and $\|w_2(p,v)\| \neq 0$, for which we define $W:=\{(p,v)\in\setC' \cap \mathcal{V}_i \mid \frac{w_1(p)}{\|w_1(p)\|}=-\frac{w_2(p,v)}{\|w_2(p,v)\|}\}$. After some algebraic manipulation, we derive the following expression from \eqref{eq:CBF_u} and \eqref{eq:CLF_u} for $(p,v) \in W$:
\begin{align}
    \text{\small $\frac{w_1(p)}{\|w_1(p)\|}u^*$} & \text{\small $\le \frac{-2v^{T}Av + \frac{\partial \alpha(h(p))}{\partial p} v + \alpha'(h'(p,v))}{\|w_1(p)\|}$} \label{eq:CBF_uu} \\
    \text{\small $\frac{w_1(p)}{\|w_1(p)\|} u^*$} & \text{\small $\ge \frac{-v^T P_2^T v - \boldsymbol{p}_i^T P_1 v  + W_3(\boldsymbol{p}_i,v)}{\|w_2(p,v)\|}.$} \label{eq:CLF_uu}
\end{align}
Then, one can define a lower-bound to $-2v^{T}Av + \frac{\partial \alpha(h(p))}{\partial p} v \ge M_1(p,v)$, and an upper-bound to $-v^T P_2^T v - \boldsymbol{p}_i^T P_1 v \le  M_2(p,v)$, for $(p,v) \in W$. The idea is now to select the function $\alpha'$ based on these upper and lower bounds so that \eqref{eq:CBF_uu} and \eqref{eq:CLF_uu} are always mutually feasible. Let us first analyze the case where  $\alpha'(h'(p,v)) = 0$, i.e. $(p,v) \in W \cap \boundC'$, where we have $ M_1(p,v) \le -2v^{T}Av - \frac{d \alpha(h(p))}{dh} \alpha(h(p)) < 0$, and existence of control inputs that satisfy the constraints can not necessarily be guaranteed anymore. To solve this, let us prove that one can design the function $\alpha$ used for the CBF such that $W \cap \boundC' = \emptyset$. If we have $(p,v) \in \boundC'$, i.e., $\alpha'(h'(p,v)) = 0$, then it follows $2(p-p_0)^TAv = \alpha(h(p))$. We also have $(p,v) \in W$, i.e. $2(p-p_0)^TA = -\frac{\|w_1\|}{\|w_2\|} (\boldsymbol{p}_i^T P_2 + v^T P_3)$. Therefore, if $(p,v) \in W \cap \boundC'$, one can infer that $(\boldsymbol{p}_i^T P_2 + v^T P_3)v = -\frac{\|w_2\|}{\|w_1\|} \alpha(h(p))$. Recall that $p \in \boundC$ are excluded from such a set, along with $p'$ in a neighborhood around $p$, so that $\alpha(h(p))>0$. Further note that since $P$ is positive definite, the quadratic form in $V(\boldsymbol{p}_i,v)$ is convex, and the gradient $\|w_2\|=0$ if and only if $\boldsymbol{p}_i=v=0$. 
Also note that $(\mathbf{x}_i,0) \not\in \boundC'$, and therefore $\|w_2\| > 0 \;\forall (p,v) \in W \cap \boundC'$.
We can now design $\alpha(h(p))$ as follows. We first lower bound $({p}_i^T P_2 + v^T P_3)v \ge M_3(p,v)$ for all $(p,v) \in \boundC'$. Therefore, we can then choose $\alpha(h(p))$ such that $M_3(p,v) + \frac{\|w_2\|}{\|w_1\|} \alpha(h(p)) \ge 0$ $\forall (p,v) \in \boundC'$ , which then implies $(\boldsymbol{p}_i^T P_2 + v^T P_3)v > -\frac{\|w_2\|}{\|w_1\|}2(p-p_0)^TAv$, and therefore $(\boldsymbol{p}_i^T P_2 + v^T P_3) \neq -\frac{\|w_2\|}{\|w_1\|} 2(p-p_0)^TA$, i.e. $w_1(p)$ and $w_2(p,v)$ are not parallel.
 
We now know that $\alpha'(h'(p,v)) \neq 0$, and even more that $h'(p,v)$ is lower bounded (on $W \cap \boundC'$). This allows one to appropriately choose $\alpha'(h'(p,v))$ and $W_3(\boldsymbol{p}_i,v)$ such that $\frac{1}{\|w_1(p)\|}\left[ M_1(p,v) + \alpha'(h'(p,v))\right]$ is sufficiently larger than $\frac{1}{\|w_2(p,v)\|} \left[M_2(p,v) + W_3(\boldsymbol{p}_i,v)\right]$ for all $(p,v) \in W$, which in turn guarantees existence of $u^* \in \Real^m$ that satisfies the constraints \eqref{eq:CBF_uu} and \eqref{eq:CLF_uu} simultaneously.
%
\qed

\emph{Proof of Theorem \ref{th:switch}:}
%
%
%
Assume that $\setC_i \cap \setC_j\neq \emptyset$ and $\setC_i \cap \setC_j$ is not a singleton. Then pick $p^*$ such that $p^*\in \text{Int}(\mathcal{C}_i)$ and $p^*\in \text{Int}(\mathcal{C}_j)$ which implies that $\alpha(h_i(p^*))>0$ and $\alpha(h_j(p^*))>0$. We now want to show that $\setC'_i \cap \setC'_j\neq \emptyset$ and $\setC'_i \cap \setC'_j$ is not a singleton.  Looking at \eqref{eq:hprime} and in order to show that $\setC'_i \cap \setC'_j\neq \emptyset$, it is enough to show that $v^*$ is the solution to the following system of equations:
\begin{align*}
    \text{\small $- 2(p^* - p_{0,i})^T A_i v^* + \alpha(h_i(p^*))$} & \text{\small $\geq 0$} \\
    \text{\small $- 2(p^* - p_{0,j})^T A_j v^* + \alpha(h_j(p^*))$} & \text{\small $\geq 0.$}
\end{align*}
If now $v^* = 0$, we see that this trivially holds. Due to continuity of $\alpha$, $h_i$, and $h_j$, it holds that there exists a $v^*\neq 0$ such that the above system of equations still holds, showing that $\setC'_i \cap \setC'_j$ is not a singleton.
\qed

\bibliographystyle{IEEEtran}
\bibliography{literature}

\end{document}